\documentclass{casc}
\newenvironment{algorithm}[1]
{
\begin{figure}[h]
  \begin{center}
   \textbf{Algorithm:} \textbf{\emph{#1}}\\*
    \begin{tabular}{p{330pt}} 
}
{
 \\ 
 \end{tabular}
 \end{center}
\end{figure}
}
\def\R{{\rm I\!R}}
\def\C{{\mathchoice {\setbox0=\hbox{$\displaystyle\rm C$}\hbox{\hbox
to0pt{\kern0.4\wd0\vrule height0.9\ht0\hss}\box0}}
{\setbox0=\hbox{$\textstyle\rm C$}\hbox{\hbox
to0pt{\kern0.4\wd0\vrule height0.9\ht0\hss}\box0}}
{\setbox0=\hbox{$\scriptstyle\rm C$}\hbox{\hbox
to0pt{\kern0.4\wd0\vrule height0.9\ht0\hss}\box0}}
{\setbox0=\hbox{$\scriptscriptstyle\rm C$}\hbox{\hbox
to0pt{\kern0.4\wd0\vrule height0.9\ht0\hss}\box0}}}}
\def\Z{{\mathchoice {\hbox{$\textstyle\sf Z\kern-0.4em Z$}}
{\hbox{$\textstyle\sf Z\kern-0.4em Z$}}
{\hbox{$\scriptstyle\sf Z\kern-0.3em Z$}}
{\hbox{$\scriptscriptstyle\sf Z\kern-0.2em Z$}}}}
\def\Q{{\mathchoice {\setbox0=\hbox{$\displaystyle\rm
Q$}\hbox{\raise
0.15\ht0\hbox to0pt{\kern0.4\wd0\vrule height0.8\ht0\hss}\box0}}
{\setbox0=\hbox{$\textstyle\rm Q$}\hbox{\raise
0.15\ht0\hbox to0pt{\kern0.4\wd0\vrule height0.8\ht0\hss}\box0}}
{\setbox0=\hbox{$\scriptstyle\rm Q$}\hbox{\raise
0.15\ht0\hbox to0pt{\kern0.4\wd0\vrule height0.7\ht0\hss}\box0}}
{\setbox0=\hbox{$\scriptscriptstyle\rm Q$}\hbox{\raise
0.15\ht0\hbox to0pt{\kern0.4\wd0\vrule height0.7\ht0\hss}\box0}}}}
\def\Fp{{\rm I\!F_{\!p}}}
\newcommand{\F}[1]{{\rm I\!F_{\!#1}}}
\def\square#1{\mathop{\mkern0.5\thinmuskip\vbox{\hrule
    \hbox{\vrule\hskip#1\vrule height#1 width 0pt\vrule}\hrule}
    \mkern0.5\thinmuskip}}
\def\Square{\mathchoice{\square{8pt}}{\square{7pt}}{\square{6pt}}{\square{5pt}}}

\newcommand{\Matrc}[5]{
\!\!\!\mbox{
\begin{tabular}{c}
\scriptsize{$#2$}\\[#5pt]
\end{tabular}
}
\hspace*{#4pt}\stackrel{\mbox{\scriptsize{$#3$}}}{#1}
}
\newcommand{\Tdim}[3]{\begin{tabular}{c}
  \\[-7pt]
  $#1$ \\
  $#2$ \\
  $#3$\\[2pt]
\end{tabular}}
\def\Tmarg{\hspace*{-1.21pt}}
\def\Tmargleftc{\hspace*{-2pt}}
\newcommand{\Ti}[1]{\Tmarg#1\Tmarg&}
\newcommand{\Tie}[1]{\Tmarg#1\Tmarg\\ \hline}
\newcommand{\Tm}[1]{$#1$&}
\newcommand{\Tc}[4]{\Tmargleftc\begin{tabular}{c}
       \\[-10pt]
  {#1} \\[-3pt]
  {\phantom{$\mathbf{#4}$} #2} $\mathbf{#4}$\\[-3pt]
  {#3} \\[-2pt]
\end{tabular}\Tmarg&}
\newcommand{\Tce}[4]{\Tmargleftc\begin{tabular}{c}
       \\[-10pt]
  {#1} \\[-3pt]
  {\phantom{$\mathbf{#4}$} #2} $\mathbf{#4}$\\[-3pt]
  {#3} \\[-2pt]
\end{tabular}\Tmarg\\}
\newcommand{\Td}[3]{\Tmarg\begin{tabular}{c}
       \\[-10pt]
  {#1} \\[-3pt]
  {#2} \\[-3pt]
  {#3} \\[-2pt]
\end{tabular}\Tmarg&}

\newcommand{\Treprow}{$\rightarrow$&}
\newcommand{\Treprowe}{$\!\rightarrow\!$\\}
\begin{document}
\title{Modular Algorithm for Computing Cohomology:\\
Lie Superalgebra of Special Vector Fields on \\
$(2|2)$-dimensional Odd-Symplectic Superspace}
\author{Vladimir V. Kornyak}
\institute{Laboratory of Information Technologies\\
           Joint Institute for Nuclear Research\\
           141980 Dubna, Russia \\
           \email{kornyak@jinr.ru}}
\maketitle
\begin{abstract}
We describe an essential improvement of our recent algorithm
for computing cohomology  of Lie (super)algebra based on partition
of the whole cochain complex into minimal subcomplexes.
We replace the arithmetic of rational numbers or integers by a much
cheaper arithmetic of a modular field and use the inequality between
the dimensions of cohomology  $H$ over any modular field $\Fp=\Z/p\Z$
and over $\Q$: $\dim H(\Fp) \geq \dim H(\Q )$. With this inequality
we can, by computing over arbitrary $\Fp,$ quickly find the (usually,
rare) subcomplexes for which $\dim H(\Fp) > 0$ and then carry out the
full computation over $\Q$ within these subcomplexes.

We also present the results of application of the corresponding
\emph{\textbf{C}} program to the Lie superalgebra of special vector
fields preserving an ``odd-symplectic'' structure on the
$(2|2)$-dimensional supermanifold. For this algebra, we found some
new basis elements of the cohomology in the trivial module.
\end{abstract}

\section{Introduction}

Recently we proposed a new algorithm for computation of cohomology
of a wide class of Lie superalgebras. This algorithm reduces the
computation for the whole cochain complex to a number of smaller
tasks within smaller subcomplexes. One can demonstrate that if $T$
is the computation time for the whole complex, then partition of the
complex into $N$ subcomplexes reduces the computation time roughly
to the value $T/N^2$. Thus, the approach appeared to be efficient
enough to cope with several difficult tasks in computing cohomology
for particular Lie (super)algebras
\cite{KornCASC01,KornProg02,KornProg02a,KornCASC02,KornProg03}.
More detailed experiments with the \emph{\textbf{C}} implementation
of the algorithm, including profiling, reveal that arithmetic operations
over $\Q$ take the main part of computation time
(usually more than 90\% for large tasks). The same is true if  $\Q$ is
replaced by $\Z$ (though computation becomes somewhat faster).

A standard way to reduce negative influence of this bottleneck is to
compute several modular images of the problem with subsequent restoring
the result over $\Z$ or $\Q$ by the Chinese remaindering or an algorithm
for recovering a rational number from its modular residue
\cite{Wang,Collins}. Though, as is clear, the sum of sizes of modules used
for constructing images can not be less than the size of maximum integer
in the final result, the modular approach allows the intermediate swelling
of coefficients to be avoided. Moreover the use of modular images is much
more advantageous in the case of (co)homology computation than in the
traditional problems of linear algebra. As we demonstrate further, the
overwhelming part of computation can be accomplished using only one
modular image. Recall that the Gauss elimination, the basic constituent of
algorithms for computing (co)homology, over $\Fp$ has only \emph{cubic}
computational complexity in contrast to the \emph{exponential} one over
$\Q$ or $\Z.$

Note that the approach presented here can be applied not only for the Lie
superalgebras, but in more general case of computation of homology or
cohomology, especially when there is a practical method of splitting
(co)chain complex into smaller subcomplexes.

To demonstrate the power of the new algorithm and program,
we present the results of computation of cohomology in the trivial module
for the algebra $\mathrm{SLe(2)}.$ This is an example (for $n=2$) of the
Lie superalgebra of special (divergence free) vector fields on the
$(n|n)$-dimensional supermanifold preserving the odd version of symplectic
(periplectic, as A.~Weil called it) structure \cite{Leites,MPI}.
These superalgebras, being super counterparts of the Lie algebras of
Hamiltonian vector fields, are vital in the Batalin--Vilkovisky formalism,
see \cite{GPS}.

\section{Combining Splitting Algorithm with Modular Search}
The $k$th cohomology is defined as the quotient group
$$H^k = Z^k/B^k \equiv {\mathrm{Ker}\ d^k}/{\mathrm{Im}\ d^{k-1}}$$
for the \emph{cochain complex}
\begin{equation}
0\to C^0\stackrel{d^0}{\longrightarrow}\cdots
\stackrel{d^{k-2}}{\longrightarrow}C^{k-1}\stackrel{d^{k-1}}
{\longrightarrow} C^k\stackrel{d^k}{\longrightarrow}C^{k+1}
\stackrel{d^{k+1}}{\longrightarrow}\cdots.
\label{whole}
\end{equation}
Here, the $C^k$ are abelian groups of cochains, graded by the integer $k$
(called \emph{dimension} or \emph{degree}); the $d^k$ are differentials
($d^k \circ d^{k-1} = 0$); the $Z^k = \mathrm{Ker}\ d^k$ and $B^k =
\mathrm{Im}\ d^{k-1}$ are the subgroups of cocycles and coboundaries,
respectively (see \cite{Fuks} for details). In order to apply without
restrictions the linear algebra algorithms, we assume that the groups of
cochains are additive groups of certain linear spaces or modules and we
shall use the corresponding terms in the subsequent text.

\subsection{Splitting Algorithm}
To compute the $k$th cohomology, it suffices to consider
the following part of (\ref{whole}):
\begin{equation}
C^{k-1}\stackrel{d^{k-1}}{\longrightarrow}
C^k\stackrel{d^k}{\longrightarrow}C^{k+1}. \label{part}
\end{equation}

First of all we split (\ref{part}) using the $\Z$-grading in the cochain
spaces induced by the gradings in the Lie (super)algebra (and the module
over this algebra) involved in the construction of the cochain spaces:
$$
\left(C^{k-1}\stackrel{d^{k-1}}{\longrightarrow}
C^k\stackrel{d^k}{\longrightarrow}C^{k+1}\right)
= \bigoplus_{g\in G} \left(C^{k-1}_g
\stackrel{d^{k-1}_g}{\longrightarrow} C^k_g
\stackrel{d^k_g}{\longrightarrow}C^{k+1}_g\right).
$$
Here, $G \subseteq \Z$ is a grading subset.

It appears that, as a rule, any subcomplex in a given degree $g$
can be split, in turn, into smaller subcomplexes:
\begin{equation}
\left(C^{k-1}_g\stackrel{d^{k-1}_g}{\longrightarrow}
C^k_g\stackrel{d^k_g}{\longrightarrow}C^{k+1}_g\right)
=\bigoplus_{s\in S} \left(C^{k-1}_{g,s}
\stackrel{d^{k-1}_{g,s}}{\longrightarrow}
C^k_{g,s}\stackrel{d^k_{g,s}}{\longrightarrow}C^{k+1}_{g,s}\right).
\label{spart}
\end{equation}
Here $S$ is a finite or infinite set of subcomplexes.

Equation (\ref{spart}) means that the spaces $C^i_g$ split
into the direct sum of subspaces
$$C^i_g = \bigoplus_{s\in S} C^i_{g,s},$$
and the matrices of the linear mappings $d^i_g$ can be represented
in the block-diagonal form
$$d^i_g = \bigoplus_{s\in S} d^i_{g,s}.
$$

The construction of these subcomplexes is
the \emph{central part} of the splitting algorithm.

Thus, the whole task reduces to a collection of easier tasks of computing
\begin{equation}
H^k_{g,s} =  {\mathrm{Ker}\ d^k_{g,s}}/{\mathrm{Im}\ d^{k-1}_{g,s}}.
\label{insubcomplex}
\end{equation}
As a basis of the cochain space $C^k_g,$ we choose the set of
super skew-symmetric monomials of the form
\begin{equation}
c(e_{i_1},\ldots,e_{i_k};a_\alpha)\equiv e^{i_1}\wedge\cdots\wedge e^{i_k}
\otimes a_\alpha. \label{monom}
\end{equation}
Here, $e_i$ and $a_\alpha$ are basis elements of the algebra and module,
respectively, and $e^i$ is the dual element to $e_i,$ that is,
$e^i(e_j) = \delta^i_j.$
The degrees of factors in (\ref{monom}) satisfy the relation
$$
\mathrm{gr}(e^{i_1})+\cdots+\mathrm{gr}(e^{i_k})
+\mathrm{gr}(a_\alpha)= g.
$$
Notice that $\mathrm{gr}(e^i)=-\mathrm{gr}(e_i)$
and this is a serious obstacle to extraction of finite-dimensional
subcomplexes for infinite-dimensional Lie (super)algebras when computing
cohomology in the adjoint module (important in the deformation theory).
We also assume that $i_1\leq\cdots\leq i_k.$

To construct a subcomplex
\begin{equation}
C^{k-1}_{g,s}\stackrel{d^{k-1}_{g,s}}{\longrightarrow}
C^k_{g,s}\stackrel{d^k_{g,s}}{\longrightarrow}C^{k+1}_{g,s}
\label{subcomplex}
\end{equation}
from the sum in right hand side of (\ref{spart}),
we begin with choosing somehow an arbitrary \emph{starting}
monomial
$m^k_{g,\mbox{\scriptsize{\,start}}}$
of the form (\ref{monom}).
There are various choices of the starting monomial and
the time and space efficiency of computation depends on
the choice. Having no better idea, we use at present the following
three strategies:
choice of a lexicographically minimal, or lexicographically maximal,
or random monomial. We call these strategies \emph{bottom},
\emph{top} and \emph{random}, respectively.
Among these, the top strategy seems to be most efficient
(see experimental data in Tables \ref{TimesH} and \ref{TimesSLe})
and it is used by default.
Nevertheless, other strategies help sometimes to push through
difficult tasks when the top strategy fails.

Then we construct the three sets
$M^{k-1}_{g,s},$ $M^k_{g,s}$ and $M^{k+1}_{g,s}$
of basis monomials for $C^{k-1}_{g,s},$ $C^k_{g,s}$ and
$C^{k+1}_{g,s},$ respectively, by the procedure
\emph{\textbf{ConstructSubcomplex}} presented on page
\pageref{ConstructSubcomplex}.
\begin{algorithm}{ConstructSubcomplex\label{ConstructSubcomplex}}
\textbf{~Input:~~}  $m^k_{g,\mbox{\scriptsize{\,start}}}$,
starting $(k,g)$-monomial\\
\textbf{~Output:} $M^{k-1}_{g,s}$, $M^{k}_{g,s}$, $M^{k+1}_{g,s}$,
monomial bases of cochain spaces in the current\\
~\phantom{\textbf{~Output:}} subcomplex $s$:
$C^{k-1}_{g,s}\stackrel{d^{k-1}_{g,s}}
{\longrightarrow}C^k_{g,s}\stackrel{d^k_{g,s}}
{\longrightarrow}C^{k+1}_{g,s},$
 such that $m^k_{g,\mbox{\scriptsize{\,start}}}
 \in M^{k}_{g,s}$\\
\textbf{~Local:~~} $\widetilde{M}^{k-1}_{g,s}\subseteq
M^{k-1}_{g,s},$
$\widetilde{M}^k_{g,s}\subseteq M^k_{g,s}$,
$\widetilde{M}^{k+1}_{g,s}\subseteq M^{k+1}_{g,s}$,
working subsets\\
~\phantom{\textbf{~Output:}}of currently ``new''
(not yet processed) monomials;\\
~\phantom{\textbf{~Output:}}$W^{k-1}_{g,s}$,
$W^k_{g,s}$, $W^{k+1}_{g,s}$, working sets of
monomials;\\
~\phantom{\textbf{~Output:}}$m^{k-1}_{g,s}$,
$m^{k}_{g,s}$, $m^{k+1}_{g,s}$, working monomials\\
\textbf{~Initial setting:}\\
~~1: $\widetilde{M}^{k-1}_{g,s}
:=M^{k-1}_{g,s}:=\emptyset$\\
~~2: $\widetilde{M}^k_{g,s}:=M^k_{g,s}
:=\{m^k_{g,\mbox{\scriptsize{\,start}}}\}$\\
~~3: $\widetilde{M}^{k+1}_{g,s}:=M^{k+1}_{g,s}:=\emptyset$\\
\textbf{~Loop over $k$-monomials:}\\
~~4:    \emph{\textbf{while}} $\widetilde{M}^k_{g,s} \neq \emptyset$
\emph{\textbf{do}}\\
~~5:~~~~ $m^k_{g,s} := \mbox{\emph{\textbf{TakeMonomialFromSet}}}
(\widetilde{M}^k_{g,s})$\\
~~~~~~~~ \textbf{Supplement the set $M^{k-1}_{g,s}$}:\\
~~6:~~~~ $W^{k-1}_{g,s}:=\mbox{\emph{\textbf{InverseImageMonomials}}}
(d^{k-1}_{g,s}, m^k_{g,s})$\\
~~7:~~~~    $M^{k-1}_{g,s}:=M^{k-1}_{g,s}\cup W^{k-1}_{g,s}$\\
~~8:~~~~    $\widetilde{M}^{k-1}_{g,s}:=\widetilde{M}^{k-1}_{g,s}
\cup W^{k-1}_{g,s}$\\
~~~~~~~~ \textbf{Supplement the set $M^{k+1}_{g,s}$}:\\
~~9:~~~~ $W^{k+1}_{g,s}:=\mbox{\emph{\textbf{ImageMonomials}}}
(d^k_{g,s}, m^k_{g,s})$\\
~10:~~~~    $M^{k+1}_{g,s}:=M^{k+1}_{g,s}\cup W^{k+1}_{g,s}$\\
~11:~~~~ $\widetilde{M}^{k+1}_{g,s}:=\widetilde{M}^{k+1}_{g,s}
\cup W^{k+1}_{g,s}$\\
~~~~~~~~    \textbf{Exclude processed monomial $m^{k}_{g,s}$}:\\
~12:~~~~ $\widetilde{M}^{k}_{g,s}:=\widetilde{M}^{k}_{g,s}
\setminus \{m^{k}_{g,s}\}$\\
~~~~~~~~ \textbf{Loop over $(k+1)$-monomials:}\\
~13:~~~~ \emph{\textbf{while}} $\widetilde{M}^{k+1}_{g,s}
\neq \emptyset$ \emph{\textbf{do}}\\
~14:~~~~~~~~  $m^{k+1}_{g,s} := \mbox{\emph{\textbf{TakeMonomialFromSet}}}
(\widetilde{M}^{k+1}_{g,s})$\\
~~~~~~~~~~~~    \textbf{Supplement the set $M^{k}_{g,s}:$}\\
~15:~~~~~~~~ $W^k_{g,s}:=\mbox{\emph{\textbf{InverseImageMonomials}}}
(d^k_{g,s}, m^{k+1}_{g,s})$\\
~16:~~~~~~~~    $M^k_{g,s}:=M^k_{g,s}\cup W^k_{g,s}$\\
~17:~~~~~~~~    $\widetilde{M}^k_{g,s}:=\widetilde{M}^k_{g,s}
\cup W^k_{g,s}$\\
~~~~~~~~~~~~    \textbf{Exclude processed monomial $m^{k+1}_{g,s}:$}\\
~18:~~~~~~~~    $\widetilde{M}^{k+1}_{g,s}:=\widetilde{M}^{k+1}_{g,s}
\setminus \{m^{k+1}_{g,s}\}$\\
~19:~~~~ \emph{\textbf{od}}\\
~~~~~~~~   \textbf{Loop over $(k-1)$-monomials:}\\
~20:~~~~ \emph{\textbf{while}} $\widetilde{M}^{k-1}_{g,s}
\neq \emptyset$ \emph{\textbf{do}}\\
~21:~~~~~~~~ $m^{k-1}_{g,s}:=\mbox{\emph{\textbf{TakeMonomialFromSet}}}
(\widetilde{M}^{k-1}_{g,s})$\\
~~~~~~~~~~~~ \textbf{Supplement the set $M^{k}_{g,s}$:}\\
~22:~~~~~~~~ $W^k_{g,s}:=\mbox{\emph{\textbf{ImageMonomials}}}
(d^{k-1}_{g,s}, m^{k-1}_{g,s})$\\
~23:~~~~~~~~ $M^k_{g,s}:=M^k_{g,s}\cup W^k_{g,s}$\\
~24:~~~~~~~~    $\widetilde{M}^k_{g,s}:=\widetilde{M}^k_{g,s}\cup W^k_{g,s}$\\
~~~~~~~~~~~~ \textbf{Exclude processed monomial $m^{k-1}_{g,s}$:}\\
~25:~~~~~~~~ $\widetilde{M}^{k-1}_{g,s}:=\widetilde{M}^{k-1}_{g,s}
\setminus \{m^{k-1}_{g,s}\}$\\
~26:~~~~ \emph{\textbf{od}}\\

~27: \emph{\textbf{od}}\\
~28: \emph{\textbf{return}} $M^{k-1}_{g,s}$, $M^{k}_{g,s}$,
$M^{k+1}_{g,s}$\\[-7pt]
\end{algorithm}

The function \emph{\textbf{TakeMonomialFromSet}}
called within \emph{\textbf{ConstructSubcomplex}} takes
the current monomial from a set of monomials.

The function \emph{\textbf{InverseImageMonomials}} generates
the set of $(q-1)$-monomials whose images with respect to
the mapping $d^{q-1}_{g,s}$ contain a given $q$-monomial.

The function \emph{\textbf{ImageMonomials}} generates
the set of $(q+1)$-monomials whose inverse images with
respect to the $d^q_{g,s}$ contain a given $q$-monomial.

In the finite-dimensional case, the loops in the procedure
\emph{\textbf{ConstructSubcomplex}} are finite and we obtain in the end
a minimal subcomplex of the form (\ref{subcomplex}). This is the unique
minimal subcomplex involving the starting monomial
$m^k_{g,\mbox{\scriptsize{\,start}}}$.

\subsection{Modular Search}

Let us consider in more detail the procedure of computation of cohomology
within the subcomplex in accordance with  formula (\ref{insubcomplex}).
From now on we assume that $C^{k-1}_{g,s},\ C^k_{g,s}$ and $C^{k+1}_{g,s}$
in (\ref{subcomplex}) are finite-dimensional spaces over $\Q$ or $\Fp$ or
modules over $\Z.$

Since important in mathematics and physics fields $\R$ and $\C$ are,
in principle, non-algorithmic objects, our main interest will be focused
on the cohomology over the field $\Q$ (or its algebraic extentions).
In accordance with a general theorem in the homological algebra,
called the \emph{universal coefficient theorem} \cite{FomenkoFuks},
(co)homology with coefficients from an arbitrary abelian group $G$
can be expressed in terms of (co)homology with coefficients in $\Z$.
Thus, we can carry out the computation over $\Z$ and then go to the
coefficient group we are interested in. Let us consider
now the connection between  $H^k_{g,s}(\Q)$ and $H^k_{g,s}(\Z).$
In the finite-dimensional case, the
group  $H^k_{g,s}(\Z)$ is a finitely generated abelian
group having the following canonical representation
\begin{equation}
H^k_{g,s}(\Z)\simeq
\overbrace{(\underbrace{\Z\oplus\cdots\oplus\Z}_{\beta^k})
}^{\mbox{\emph{free part}}}\oplus
\overbrace{\Z_{t_1}\oplus\cdots\oplus\Z_{t_r}}^{\mbox{\emph{torsion}}}.
    \label{abelgroup}
\end{equation}
Here,  $\beta^k$, the number of copies of the integer group $\Z$,
is called the \emph{rank} of the abelian group $H^k_{g,s}(\Z)$ or
the \emph{Betti number}. The cyclic groups $\Z_{t_i}$ are called
the \emph{torsion subgroups} and their orders  $t_i$,
having the property  $t_i > 1,$ $t_1|t_2,$ $t_2|t_3,\ldots$ and so on,
are called the \emph{torsion coefficients}.

In the case of cohomology, the universal coefficient theorem is
expressed by the following split short exact sequence
\begin{equation}
0 \rightarrow H^k_{g,s}(\Z)\otimes G \rightarrow
H^k_{g,s}(G) \rightarrow \mathrm{Tor}(H^{k+1}_{g,s}(\Z),G)
\rightarrow 0,
\label{exact}
\end{equation}
where the operation $\mathrm{Tor}$ is the
\emph{periodic product} of abelian groups.
In our context
$$
\mathrm{Tor}(H^{k+1}_{g,s}(\Z),G) =
(\mbox{\emph{torsion}}\ H^{k+1}_{g,s}(\Z))\otimes(\mbox{\emph{torsion}}\ G).
$$
The term ``split'', in application to sequence (\ref{exact}),
means the possibility to construct the isomorphism
\begin{equation}
H^k_{g,s}(G) \simeq H^k_{g,s}(\Z)\otimes G
\oplus \mathrm{Tor}(H^{k+1}_{g,s}(\Z),G).
\label{isomorphism}
\end{equation}
Replacing $G$ by $\Q$ in (\ref{isomorphism})
and taking into account that $\mathrm{Tor}(A,\Q)= 0$
for any abelian group $A$, we have
$H^k_{g,s}(\Q) \simeq H^k_{g,s}(\Z)\otimes \Q$.
Since $\Z_m\otimes\Q=0$ for arbitrary $m$
and $\Z\otimes\Q\simeq\Q,$ the dimension of
$H^k_{g,s}(\Q),$ interpreted as vector space over
$\Q$, coincides with the rank (Betti number) $\beta^k$
of the group $H^k_{g,s}(\Z).$

Our modular approach is based on the following
\begin{theorem}
\label{theorem}
\begin{equation}
\dim H^k_{g,s}(\Fp) \geq \dim H^k_{g,s}(\Q). \label{inequality}
\end{equation}
\end{theorem}
\textbf{Remarks:}
\begin{enumerate}
    \item Inequality  (\ref{inequality}) means that non-trivial
    cohomology classes computed over the field of rational numbers $\Q$
    can exist only in the subcomplexes with non-trivial cohomology classes
    computed over the finite field $\Fp$ with arbitrary prime $p.$
    \item H. Khudaverdian turned author's
    attention to the fact that inequality (\ref{inequality}) can be deduced
    immediately from the universal coefficient theorem:
    considering the product $H^k_{g,s}(\Z)\otimes \Fp$ and taking into account
    representation (\ref{abelgroup}) and isomorphism $\Z\otimes\Fp\simeq\Fp$,
    we see that the dimension of $H^k_{g,s}(\Fp)$, as a vector space
    over $\Fp$, can not be less than $\beta^k$ (only additional dimensions
    may appear, if the torsions in $H^k_{g,s}(\Z)$ or $H^{k+1}_{g,s}(\Z)$
    contain cyclic groups of the form $\Z_{p^m}$).

  Nevertheless, we give here a direct constructive proof in order to
  demonstrate in parallel the main ideas of (co)homology computation.
\end{enumerate}
\emph{Proof.}
To prove inequality (\ref{inequality}), we have to compute
(\ref{insubcomplex}) in such a way as to avoid cancellations of integers
and apply the modular homomorphism $\phi_p: \Z \rightarrow \Fp $
at the end of computation. Thus, it is convenient to consider
(\ref{insubcomplex}) over $\Z$ instead of $\Q.$

We assume that $p$ is odd and use a symmetric representation
of $\Fp$, i. e.,
$$\Fp=\left\{-\frac{p-1}{2},
\ldots,-1,0,1,\ldots,
\frac{p-1}{2}\right\}.
$$
We will also apply $\phi_p$ component-wise to multicomponent objects
over $\Z,$ like vectors and matrices.
\\
We begin with the following setup:
\begin{itemize}
    \item
 $\ C^{k-1}_{g,s},\ C^k_{g,s}$ and $\ C^{k+1}_{g,s}$ are represented
 as finite-dimensional modules
 $\ M^{-}=\Z^n,\ M=\Z^m,\ M^{+}=\Z^l,$ respectively, i. e.,
 $\dim C^{k-1}_{g,s} = n,\ \dim C^k_{g,s} = m,\
 \dim C^{k+1}_{g,s} = l.$\\[-6pt]
    \item the differentials $d^{k-1}_{g,s}$ and $d^k_{g,s}$ are
    represented (in the monomial bases of the form (\ref{monom}), in our case)
    as integer $m\times n$  and  $l\times m$ matrices
$$
D' \equiv  \Matrc{D'}{m}{n~}{-7}{4}= \left[
\begin{tabular}{ccc}
    $\left(d'\right)^1_1$ &  $\cdots$ & $\left(d'\right)^1_n$\\
  $\cdots$ &  $\cdots$ & $\cdots$\\
  $\left(d'\right)^m_1$ &  $\cdots$ & $\left(d'\right)^m_n$
\end{tabular}
\right]
\   \mbox{and} \
D \equiv  \Matrc{D}{l}{m}{-7}{3}= \left[
\begin{tabular}{ccc}
    $d^1_1$ & $\cdots$ & $d^1_m$\\
  $\cdots$ &  $\cdots$ & $\cdots$\\
  $d^l_1$ &  $\cdots$ & $d^l_m$
\end{tabular}
\right],
$$ respectively. We write $\Matrc{A}{i}{j}{-7}{2}$  to indicate
that matrix $A$ has $i$ rows and $j$ columns.
 \item the matrices $D$ and $D'$ satisfy the relation $DD' = 0.$
\end{itemize}
The computation of cohomology, i.e., construction of quotient module,
can be reduced to the construction of so-called
\emph{(co)homology decomposition} \cite{HiltonWylie} based on
the computation  of the Smith normal forms \cite{Birkhoff} of
the matrices representing differentials.

First of all let us determine the cocycle submodule, i.e.,
$\mathrm{Ker}\ D\subseteq M,$  by reducing the matrix $D$ to
the \emph{integer} Smith normal form $S = UDV.$ The matrix $S$
has the form
\begin{equation}
S =  \Matrc{S}{l}{~m}{-8}{2} = \left[
        \begin{tabular}{ll}
        $~~\Matrc{\widetilde{S}}{r}{~r}{-7}{2.5}$
        &
        $~~~\Matrc{O}{r}{m\!\!-\!r}{-9}{2}$
        \\
        $\Matrc{O}{l\!\!-\!r}{~r}{-6}{1.5}$
        &
        $\ \Matrc{O}{l\!\!-\!r}{m\!\!-\!r}{-10}{1.5}$
        \end{tabular}
\right], \label{s0}
\end{equation}
where, $\Matrc{O}{i}{j}{-7}{2}$  is the $i\times j$ zero matrix,
$r = \mathrm{rank}_\Z D,$ $\Matrc{\widetilde{S}}{r}{~r}{-7}{2.5}
=  \mathrm{diag}\left(s_1,\ \ldots\ ,\ s_r\right),$
$s_1,\ldots,s_r$ are positive integers
called the \emph{invariant factors} of $D.$ These invariant factors
have the property $s_i|s_{i+1}$ for all $i$. Note that there is connection
between the invariant factors and the torsion coefficients from formula
(\ref{abelgroup}), namely, the prime divisors of the invariant factors
are also divisors of some torsion coefficients. The transformation matrices
$U = \Matrc{U}{l}{l}{-7}{1.5}$ and $V = \Matrc{V}{m}{m}{-7}{1.5}$
are \emph{unimodular} integer matrices, i.e.,
$\det U = \pm 1,\ \det V = \pm 1.$ With such determinants,
these matrices are invertible and their inverses $U^{-1}$ and $V^{-1}$
are obviously integer matrices too.

Now we should consider the coboundary submodule
$$
\mathrm{Im}\ D'\subseteq \mathrm{Ker}\ D\subseteq M.
$$
Combining the relation
$$
SV^{-1}D' = UDD' = 0
$$
with the structure of the matrix $S$ (see formula (\ref{s0})),
we can reduce the matrix $D'$ determining coboundaries to the matrix
$\widetilde{D'}$ acting in the submodule of cocycles:
$$
V^{-1}D' =
\left[
        \begin{tabular}{l}
        $~~~\Matrc{O}{r}{n}{-6}{2}$
        \\
        $\Matrc{\widetilde{D'}}{m\!\!-\!r}{n}{-7}{4}$
        \end{tabular}
\right].
$$
Computing the Smith normal form $S' = \widetilde{U'}\widetilde{D'}V'$
for the reduced coboundary matrix we get
$$
S' = \Matrc{S'}{m\!\!-\!r}{n}{-7}{2.5} = \left[
        \begin{tabular}{ll}
        $~~~~~~\Matrc{\widetilde{S'}}{r'}{~r'}{-8}{2}$
        &
        $~~~~~~~\Matrc{O}{r'}{~n\!\!-\!r'}{-11}{2}$\\
        $\Matrc{O}{m\!\!-\!\!r\!\!-\!\!r'}{\,~r'}{-9}{2}$
        &
    $~\Matrc{O}{m\!\!-\!\!r\!\!-\!\!r'}{~n\!\!-\!r'}{-12}{2}$
        \end{tabular}
\right],
$$
where $r' = \mathrm{rank}_\Z \widetilde{D'}$ and
$\Matrc{\widetilde{S'}}{r'}{~r'}{-8}{2.5}  =
\mathrm{diag}\left(s'_1,\ \ldots\ ,\ s'_{r'}\right).$

We can extend the transformation matrix
$\widetilde{U'} = \Matrc{\widetilde{U'}}{m\!-\!r}{m\!-\!r}{-10}{3}$
acting in the submodule of cocycles to the transformation matrix acting
in the whole module $M:$
$$
U' = \Matrc{U'}{m}{m~}{-7}{4} = \left[
        \begin{tabular}{ll}
        $~~~~~\Matrc{I}{r}{\,r}{-7}{2.5}$
        &
        $~~~~~\Matrc{O}{r}{m\!-\!r}{-9}{2.5}$
        \\
        $\Matrc{O}{m\!-\!r}{\,r}{-6}{2.5}$
        &
    $~\Matrc{\widetilde{U'}}{m\!-\!r}{m\!-\!r}{-10}{2.5}$
        \end{tabular}
\right].
$$
Here, $\Matrc{I}{r}{\,r}{-8}{2.5}$ is the $r\times r$ identity matrix.
Using the transformation matrices $U'$ and $V$ we can transform
the initial (monomial in our case) basis $\mathbf{e} = (e_1,\ldots, e_m)$
in the module $M$ into the basis
$\mathbf{a} = (a_1,\ldots, a_m) = \mathbf{e}V\left(U'\right)^{-1}$
making the cohomology decomposition explicit
\begin{equation}
M = \underbrace{\left\langle a_1,\ldots,
    a_r\right\rangle}_{non-cocycles}\oplus
    \underbrace{\overbrace{\left\langle a_{r+1},\ldots,
    a_{r+r'} \right\rangle}^{coboundaries}
    \oplus\overbrace{\left\langle a_{r+r'+1},\ldots,
    a_m \right\rangle}^{cohomology}}_{cocycles}.
    \label{decomposition}
\end{equation}
In this decomposition we have
$$
\mathrm{Ker}\ d^k_{g,s} = \left\langle a_{r+1},\ldots,
a_m\right\rangle
$$
and
$$
\mathrm{Im}\ d^{k-1}_{g,s} = \left\langle a_{r+1},\ldots,
a_{r+r'}\right\rangle.
$$
The formula for the dimension of cohomology (Betti number)
follows from decomposition  (\ref{decomposition})
\begin{equation}
\dim H^k_{g,s}(\Q) = \beta^k = m - r - r'. \label{dimension}
\end{equation}

Now let us consider how (\ref{dimension}) changes under $\phi_p.$
The image of (\ref{dimension}) takes the form
\begin{equation}
\dim H^k_{g,s}(\Fp) = m_p - r_p - r'_p.
\label{dimensionp}
\end{equation}
Since $\phi_p$ is a ring homomorphism, we have for arbitrary
unimodular matrix $A$ with integer entries
$$
\det \phi_p(A) = \phi_p(\det A) = \phi_p(\pm1)= \pm1,
$$
that is $\phi_p$ maps the above transformation matrices into
invertible matrices. Hence the number of elements in the
decomposition basis $\mathbf{a}$ remains unchanged, $m_p = m.$
On the other hand, the invariant factors $s'_1,\ldots,s'_{r'}$
and $s_1,\ldots,s_{r}$ of the matrices $S'$ and  $S$ divisible
by $p$ vanish, hence
$\ \mathrm{rank}\ \phi_p(S) = r_p \leq r\ \mbox{ and }\
\mathrm{rank}\ \phi_p(S') = r'_p \leq r'$ and inequality (\ref{inequality})
is proved by comparing (\ref{dimension}) and (\ref{dimensionp}).
\hfill$\Square$

\subsection{Implementation}

An algorithm based on the above ideas was implemented in
the \emph{\textbf{C}} language. The program called
\emph{\textbf{LieCohomologyModular}} has the following structure:
\begin{enumerate}
    \item Input Lie (super)algebra $A$, module $X$ over $A$,
    cohomology degree (dimension) $k$ and grade $g$.
    $A$ and $X$ should be defined over (some algebraic extension of)
    $\Z$ or $\Q.$
    \item Construct the full set $M^k_g$ of $k$-cochain monomials in grade $g$.
    \item Choose a prime $p$ for searching subcomplexes with non-trivial
    cocycles by computing over $\Fp$.
    \item \label{choosemonom} Choose an element $m^k_g \in M^k_g$
    (the starting monomial).
    \item   Construct a minimal subcomplex $s:
C^{k-1}_{g,s}\stackrel{d^{k-1}_{g,s}}{\longrightarrow}
C^k_{g,s}\stackrel{d^k_{g,s}}{\longrightarrow}
C^{k+1}_{g,s}$ such that $m^k_g \in C^k_{g,s}.$
    \item  Compute $n = \dim H^k_{g,s}(\Fp)$.
    \item \label{computeh} If $n > 0,$  then compute $H^k_{g,s}$ over
    $\Z$ or $\Q$ (or their extensions).
  We can use here the Chinese remaindering or the rational recovery algorithm
  as more efficient procedures than direct computation over $\Z$ or $\Q$.
    \item  Delete all basis monomials of $C^k_{g,s}$ from  $M^k_g.$
    \item  If $M^k_g$ is empty, then stop computation, otherwise go to
    Step \ref{choosemonom}.
\end{enumerate}
In the current implementation we obtain the relations determining cocycles
and coboundaries (in fact, the rows of matrices of differentials) within
the procedure \emph{\textbf{ConstructSubcomplex}}. These relations are
generated one by one as by-product of the functions
\emph{\textbf{InverseImageMonomials}} and  \emph{\textbf{ImageMonomials}}.
To prevent unnecessary memory consumption, every newly arising relation is
reduced modulo the system of relations existing to the moment and, if
the result is not zero, the new relation is added to the system. Thus, we
automatically have the matrices of differentials in the normal form just
after completion of the procedure \emph{\textbf{ConstructSubcomplex}}.
This process is obviously equivalent to the Gauss elimination method, the
most standard method for the computation of the Smith normal form of a
matrix.

In recent years, a number of new fast algorithms for the determination
of Smith normal form have been elaborated \cite{DSV01,DHSW03}.
These algorithms appear to be well suited to the (co)homology computation.
It is worth to study the possibility to incorporate these algorithms in our
implementation. We could, for example, remove the generation of relations from
the functions
\emph{\textbf{InverseImageMonomials}} and \emph{\textbf{ImageMonomials}}
making them as fast as possible. Then, after construction of subcomplex
with the help of these modified functions, we should generate the matrices
of differentials separately and apply the fast algorithms to these matrices.
Of course, this modification should be done if the total computation time
decreases without substantial increase in the memory consumption.
The works \cite{DSV01,DHSW03} contain a detailed analysis of the properties
of the sets of primes most appropriate for application of modular algorithms
to a given matrix.

Here we give only a few comments concerning the choice of prime $p$
in our algorithm.
These comments are based mainly on experiments with the program.

From the practical point of view, we should use only primes $p$ for which
all operations in $\Fp$ can be done within one machine word. Thus, for
32bit architecture we should choose $p$ from the set of 8951 primes
$(3,5,\ldots,92681).$ A good choice should not produce excessive cocycles.
Of course, such cocycles will be removed at Step \ref{computeh} anyway,
but at the expense of additional work. In our context, an ``unlucky''
prime is the one which divides the invariant factors of matrices of
differentials (or, in other words, the torsion coefficients of cohomology
over $\Z$) and, as is clear, the probability for a given prime to be
unlucky diminishes as the prime grows. On the other hand, there is an
increase of time (in the examples we have computed, up to factor 2 or 3)
and space expenditures with increase of $p$ within the set
$(3,5,\ldots,92681),$ so it makes sense not to use too large primes for
searching subcomplexes with potentially non-trivial cohomology. In our
practice, we use, as a rule, a compromise: a prime near the half of 32bit
word, namely, $p = 65537 = 2^{2^4}+1,$ i. e., the 4th Fermat number.

However, quite satisfactory results can be obtained even with much
smaller primes, as is illustrated in Table
\ref{tab:dimH}. The symbols $n_{p(\Q)}$ in the boxes of this table
mean that (non-zero) $n = \dim_{\Fp(\Q)}H(\mathrm{H(2)})^k_g$,
whereas $H(\mathrm{H(2)})^k_g$ means the cohomology with coefficients
in the trivial module for the Lie algebra $\mathrm{H(2)}$ of
Hamiltonian vector fields on the $2$-dimensional symplectic manifold.
We performed computation over all modular fields $\F{3}$ through $\F{17}.$
As is seen in Table \ref{tab:dimH},
the results for $\F{17}$ fully coincide with those for $\Q$ for
all computed grades $g \in [-2,\ldots,8]$
(and for all cohomology degrees $k$). The table also illustrates
Theorem \ref{theorem}, i.e., all boxes containing
non-zero dimensions for the field $\Q$ contain also non-zero
($\geq$ same for $\Q$) dimensions for all fields $\Fp$ considered.

In Tables \ref{TimesH} and \ref{TimesSLe} we present
(considering both algebra and  superalgebra cases)
the running times for computation over $\F{17}$ of cohomology
$H^k_g(\mathrm{H(2)})$ (for $k = 7,\ 4 \leq g \leq 8$)
and  $H^k_g(\mathrm{SLe(2)})$ (for $k = 6,\ 0 \leq g \leq 4$).
The columns presented in these tables are:
the dimensions of cochain spaces, i. e.,
the sizes of matrices of differentials;
the running times in seconds for the \emph{top} strategy
of the choice of the starting monomial
and comparison of the \emph{bottom} and \emph{top} strategies.

The times in both tables were obtained on a 1133MHz Pentium III PC with 512Mb.
Note that the maximum memory consumption is near 46Mb and near 14Mb for the
tasks in Table \ref{TimesH} and in Table \ref{TimesSLe}, respectively.

\begin{table}[h!]
    \caption{Timing for $H^k_g(\mathrm{H(2)},\F{17})$, $k=7$}
    \label{TimesH}
        \begin{center}
        \begin{tabular}{c|c|c|c|c|c|}
$~g~$ & $~\dim C^{k-1}_g~$ & $~\dim C^{k}_g~$ & $~\dim C^{k+1}_g~$ & $~T_{top}~$ & $~T_{bottom}/T_{top}~$
\\ \hline
4 & ~1580 & ~1128 & ~~479 & $<1$ & 2.0
\\ \hline
5 & ~3382 & ~2730 & ~1388 & ~~~4 & 2.8
\\ \hline
6 & ~6734 & ~6132 & ~3606 & ~~27 & 3.3
\\ \hline
7 & 12766 & 12818 & ~8546 & ~214 & 3.7
\\ \hline
8 & 23074 & 25488 & 18963 & 1128 &  4.5
\\ \hline
        \end{tabular}
    \end{center}
\end{table}

\begin{table}[h!]
    \caption{Timing for $H^k_g(\mathrm{SLe(2)},\F{17})$, $k=6$}
    \label{TimesSLe}
        \begin{center}
        \begin{tabular}{c|c|c|c|c|c|}
$~g~$ & $~\dim C^{k-1}_g~$ & $~\dim C^{k}_g~$ & $~\dim C^{k+1}_g~$ & $~T_{top}~$ & $~T_{bottom}/T_{top}~$
\\ \hline
0 & ~1867 & ~6605 & ~22119 & ~~1 & 2.2
\\ \hline
1 & ~3528 & 12162 & ~39796 & ~~4 & 4.1
\\ \hline
2 & ~6546 & 22102 & ~70817 & ~21 &  4.7
\\ \hline
3 & 11878 & 39652 & 124768 & ~87 &  5.7
\\ \hline
4 & 21073 & 70110 & 217696 & 413 & 6.0
\\ \hline
        \end{tabular}
    \end{center}
\end{table}
\begin{table}[h]
    \caption{$\dim H^k_g(\mathrm{H(2)},R)$ for $R = \Q$ and $\Fp,\ p = 3,5,7,11,13,17;$
     $(k,g) \in [1,\ldots,\infty)\times[-2,\ldots,8]$}
    \label{tab:dimH}
    \begin{center}
        \begin{tabular}{l|c|c|c|c|c|c|c|c|c|c|c}
        $g\backslash k$
  & ~1~ & ~~2~ & 3 & 4 & 5 & 6 & 7 & 8 & 9 & 10 & ~11 \\ \hline
-2&   &\Tdim{1_{3}1_{5}1_{7}}{1_{11}1_{13}}{\mathbf{\{1_{17}1_{\Q}\}}}
          &   &   &\Tdim{1_{3}1_{5}1_{7}}{1_{11}1_{13}}{\mathbf{\{1_{17}1_{\Q}\}}}
                      &   &   &   &   &    &    \\ \hline
-1&   &\Tdim{}{2_{3}}{}
          &\Tdim{}{2_{3}}{}
              &\Tdim{}{2_{5}}{}
                  &\Tdim{}{2_{3}2_{5}}{}
                      &\Tdim{}{2_{3}}{}
                          &   &   &   &    &    \\ \hline
~0&   &   &\Tdim{}{1_{3}}{}
              &\Tdim{}{1_{3}}{}
                  &\Tdim{}{3_{3}}{}
                      &\Tdim{}{3_{3}}{}
                          &\Tdim{1_{3}1_{5}1_{7}}{1_{11}1_{13}}{\mathbf{\{1_{17}1_{\Q}\}}}
                              &   &   &    &    \\ \hline
~1&   &\Tdim{}{2_{5}}{}
          &\Tdim{}{2_{5}}{}
              &\Tdim{}{2_{3}2_{7}}{}
                  &\Tdim{}{8_{3}2_{7}}{}
                      &\Tdim{}{6_{3}}{}
                          &   &   &   &    &    \\ \hline
~2&\Tdim{}{~1_{3}}{}
      &\Tdim{}{4_{3}}{}
          &\Tdim{}{3_{3}}{}
              &\Tdim{}{1_{3}3_{5}}{}
                  &\Tdim{}{8_{3}4_{5}}{}
                      &\Tdim{}{11_{3}1_{5}}{}
                          &\Tdim{}{4_{3}}{}
                              &\Tdim{}{1_{3}}{}
                                  &\Tdim{}{1_{3}}{}
                                      &    &    \\ \hline
~3&   &\Tdim{}{2_{7}}{}
          &\Tdim{}{2_{3}2_{7}}{}
              &\Tdim{}{6_{3}2_{5}}{}
                  &\Tdim{}{10_{3}2_{5}}{}
                      &\Tdim{}{10_{3}}{}
                          &\Tdim{}{8_{3}2_{5}}{}
                              &\Tdim{}{6_{3}2_{5}}{}
                                  &\Tdim{}{2_{3}}{}
                                      &    &    \\ \hline
~4&   &   &\Tdim{}{5_{3}1_{5}}{}
              &\Tdim{}{20_{3}4_{5}}{}
                  &\Tdim{}{31_{3}3_{5}}{}
                      &\Tdim{}{17_{3}1_{5}3_{7}}{}
                          &\Tdim{}{2_{3}4_{5}4_{7}}{}
                              &\Tdim{}{4_{3}3_{5}1_{7}}{}
                                  &\Tdim{}{3_{3}}{}
                                      &    &    \\ \hline
~5&\Tdim{}{~2_{3}}{}
      &\Tdim{}{6_{3}}{}
          &\Tdim{}{4_{3}}{}
              &\Tdim{}{8_{3}2_{7}}{2_{11}}
                  &\Tdim{}{34_{3}2_{7}}{2_{11}}
                      &\Tdim{}{42_{3}2_{5}}{}
                          &\Tdim{}{20_{3}6_{5}2_{7}}{2_{11}}
                              &\Tdim{}{10_{3}4_{5}2_{7}}{2_{11}}
                                  &\Tdim{}{6_{3}}{}
                                      &    &    \\ \hline
~6&\Tdim{}{~1_{5}}{}
      &\Tdim{}{4_{5}}{}
          &\Tdim{}{3_{3}3_{5}}{}
              &\Tdim{}{8_{3}3_{7}}{}
                  &\Tdim{}{19_{3}4_{7}}{}
                      &\Tdim{}{35_{3}1_{7}}{}
                          &\Tdim{}{45_{3}3_{5}}{}
                              &\Tdim{}{38_{3}4_{5}}{}
                                  &\Tdim{}{17_{3}1_{5}}{}
                                      &\Tdim{}{3_{3}}{}
                                           &    \\ \hline
~7&   &\Tdim{}{2_{3}2_{11}}{}
          &\Tdim{}{16_{3}2_{5}}{2_{11}}
              &\Tdim{}{48_{3}10_{5}}{2_{13}}
                  &\Tdim{}{64_{3}14_{5}}{2_{13}}
                      &\Tdim{}{36_{3}6_{5}2_{7}}{}
                          &\Tdim{}{30_{3}2_{5}2_{7}}{2_{13}}
                              &\Tdim{}{60_{3}2_{5}}{2_{13}}
                                  &\Tdim{}{44_{3}2_{7}}{2_{11}}
                                      &\Tdim{}{8_{3}2_{7}}{2_{11}}
                                           &    \\ \hline
~8&\Tdim{}{~3_{3}}{}
      &\Tdim{}{8_{3}}{}
          &\Tdim{}{8_{3}1_{5}1_{7}}{}
              &\Tdim{}{27_{3}4_{5}4_{7}}{}
                  &\Tdim{}{87_{3}3_{5}3_{7}}{}
                      &\Tdim{}{110_{3}5_{5}1_{7}}{}
                          &\Tdim{73_{3}18_{5}5_{7}}{1_{11}1_{13}}{\mathbf{\{1_{17}1_{\Q}\}}}
                              &\Tdim{}{65_{3}15_{5}3_{7}}{}
                                  &\Tdim{}{56_{3}3_{5}}{}
                                      &\Tdim{18_{3}1_{5}1_{7}}{1_{11}1_{13}}{\mathbf{\{1_{17}1_{\Q}\}}}
                                           &\Tdim{}{1_{3}}{} \\ \hline
        \end{tabular}
    \end{center}

\end{table}
\section{Computing $H^k_g(\mathrm{SLe(2)})$}
In this section we present the results of application of the
program \textit{\textbf{LieCohomologyModular}}
to the Lie superalgebra of special vector fields preserving
periplectic structure on $(2|2)$-dimensional superspace. Periplectic
supermanifolds with a fixed volume element play an important role
in the geometrical formulation of the Batalin-Vilkovisky
formalism \cite{BV1,BV2},
an efficient method for quantizing gauge theories.

Recall that a \textit{periplectic} or an
\textit{odd symplectic manifold} is an $(n|n)$-dimensional
supermanifold equipped with an odd symplectic structure,
that is an odd non-degenerate closed 2-form. In an
analog of Darboux coordinates \cite{Shander},  it
takes the shape
\begin{equation}
\omega = \sum_{i=1}^n dx^i\wedge d\theta_i.
\label{2form}
\end{equation}
Here, $x^1,\ldots,x^n$ and $\theta_1,\ldots,\theta_n$ are even
and odd (Grassmann) variables, respectively.
The vector fields preserving 2-form (\ref{2form}) form a Lie
superalgebra denoted by $\mathrm{Le(n)}.$ The elements of
$\mathrm{Le(n)}$ can be expressed in terms of \textit{generating functions}
(also called \textit{hamiltonians}) and these generating functions generate
a nontrivial central extension of $\mathrm{Le(n)}$
called the \textit{Buttin algebra} and denoted by $\mathrm{B(n)}$.
The bracket for arbitrary two hamiltonians $f$ and $g$ is called the
\textit{Buttin bracket} or \textit{antibracket} or
\textit{odd Poisson bracket} and takes the form
\begin{equation}
 \{f,g\} =\sum_{i=1}^{n}\left(\frac{\partial f}{\partial x^i}
 \frac{\partial g}{\partial \theta_i}
 +(-1)^{p(f)}\frac{\partial f}{\partial \theta_i}
 \frac{\partial g}{\partial x^i}\right).
\end{equation}
Here, $p(f)$ is parity of function $f.$

The odd symplectic structure is a super version of the ordinary
symplectic structure
$\omega = \sum_{i=1}^n dq^i\wedge dp_i,$ where $q^i$ and $p_i$
are both even.
In symplectic case there is an invariant volume form
$\rho_\omega = \omega^n$
(Liouville theorem) and all the vector fields
preserving $\omega$ are automatically divergence free.
Contrariwise, in periplectic
case the volume is not preserved (see \cite{MPI}, for more
geometric consideration of the subject see
\cite{Ovik}) and one can impose the divergence-free
condition additionally.

Thus, we come to the \textit{special Buttin algebra} $\mathrm{SB(n)}.$ Its generating functions $f$ satisfy
the divergence free condition
\begin{equation}
\Delta f = 0,
\label{Df}
\end{equation}
 where $\Delta$ is the \textit{odd Laplacian}
\begin{equation}
\Delta = \sum_{i=1}^{n}\frac{\partial^2}{\partial x^i\partial \theta_i}.
\label{Delta}
\end{equation}
This $\Delta$-operator is, actually, the Fourier
transform with respect to the odd variables of the
usual de Rham differential (see \cite{MPI}).  Together
with its homology, $\Delta$ plays the key role in the formulation
of so-called  Batalin-Vilkovisky ``master equation''.
We can slightly reduce the special Buttin algebra by removing
constants from generating functions, i.e., taking the quotient
algebra modulo the center $\mathcal{Z}$.
The resulting algebra $\mathrm{SLe(n)}=\mathrm{SB(n)}/\mathcal{Z}$
is called the \textit{special Leites algebra.}

Since the above algebras are infinite-dimensional, in order to
compute cohomology, we introduce a grading
by prescribing grades to the variables $\{x^i\}$ and $\{\theta_i\}$
with subsequent extension of the grading to the polynomial functions
of these variables. Most natural grading can be provided by setting
$\mathrm{gr}(x^i) = 1$ and  $\mathrm{gr}(\theta_i) = -1.$ Here the
word ``natural'' means that, with this grading, the algebra  $\mathrm{B(n)}$
(and $\mathrm{Le(n)}=\mathrm{B(n)}/\mathcal{Z}$) contains
the \textit{inner grading element}
\begin{equation}
     \sum_{i=1}^{n}x^i\theta_i,
     \label{G}
\end{equation}
guaranteeing that all non-trivial cohomology classes lie in the
zero grade cochain subspaces (see \cite{Fuks}). Unfortunately,
there is no good inner grading element\footnote{
In a private communication I. Shchepochkina suggested
a general view on the inner grading elements for the algebras
with antibracket.
She noticed that $\sum_{i=1}^{n}g_i x^i\theta_i$ can be used
as the inner grading element in $\mathrm{B(n)}$ or $\mathrm{Le(n)}$
as well as (\ref{G}) assuming that $g_i$ are arbitrary integers and
$\mathrm{gr}(x^i) = g_i$ and  $\mathrm{gr}(\theta_i) = -g_i.$
This inner grading element can belong to a divergence-free
subalgebras $\mathrm{SB(n)}$ or $\mathrm{SLe(n)}$ only if
$\sum_{i=1}^{n}g_i=0$ (in virtue of equation (\ref{Df})), i. e.,
the set $\{g_i\}$ must contain non-zero
integers with opposite signs. This leads to impossibility to construct
finite-dimensional cochain subspaces of a given degree.}
in the divergence free algebra $\mathrm{SB(n)}$ (and $\mathrm{SLe(n)}$).
That is why the computation of cohomology for these algebras is much more
difficult task than for the algebras without divergence free condition.

Let us now turn to $\mathrm{SLe(2)}.$ Its  basis elements up to grade 1 are
\[
\begin{array}{cl}
\mbox{~~~~\emph{grades}~~~~}
   & \mbox{\emph{basis elements}}\\[3pt]
-2 & \left\{
\begin{array}{ccc}
    O_1 & = & \theta\psi, \\
\end{array} \right.
 \\[1pt]
-1 &
\left\{
\begin{array}{ccc}
    E_2 & = & \theta, \\
    E_3 & = & \psi,
\end{array} \right.
\\[7pt]
 0 &
\left\{
\begin{array}{ccl}
E_{4}&=&y\theta,\\
E_{5}&=&y\psi-x\theta,\\
E_{6}&=&x\psi,\\
\end{array} \right.
\\[13pt]
 1 &
\left\{
\begin{array}{ccl}
O_{7}&=&y,\\
O_{8}&=&x,\\
E_{9}&=&y^2\theta,\\
E_{10}&=&y^2\psi-2xy\theta,\\
E_{11}&=&xy\psi-\frac{1}{2}x^2\theta,\\
E_{12}&=&x^2\psi,\\
\end{array} \right.
\\[16pt]
 2 &
\left\{
\begin{array}{ccl}
\phantom{E_{12}} & ~\vdots &\\
\end{array} \right.
\end{array}
\]
Here, $x, y$ and $\theta, \psi$ are even and odd variables
of the $(2|2)$-dimensional superspace, respectively;
$E_i$ and $O_i$ are even and odd basis elements of the
Lie superalgebra, respectively. Notice that odd generating
function corresponds to the even element of superalgebra and
vice versa. This is a property of Lie superalgebras
with antibrackets called \emph{parity shift}.

Let us present also the initial part of the multiplication
table in $\mathrm{SLe(2)}$
(given are only non-zero brackets)
\begin{center}
\begin{tabular}{lcl}
$[E_{2}, E_{5}]$ & $=$ & $E_{2},$\\
$[E_{2}, E_{6}]$ & $=$ & $-E_{3},$\\
$[E_{3},E_{4}]$&$=$&$-E_{2},$\\
$[E_{3},E_{5}]$&$=$&$-E_{3},$\\[5pt]
$[E_{5},E_{4}]$&$=$&$-2E_{4},$\\
$[E_{6},E_{4}]$&$=$&$E_{5},$\\
$[E_{5},E_{6}]$&$=$&$2E_{6},$\\[5pt]
$[O_{1},O_{7}]$&$=$&$-E_{2},$\\
$[E_{5},O_{7}]$&$=$&$-O_{7},$\\
$[E_{6},O_{7}]$&$=$&$-O_{8},$\\
$[O_{1},O_{8}]$&$=$&$E_{3},$\\
$[E_{4},O_{8}]$&$=$&$-O_{7},$\\
$[E_{5},O_{8}]$&$=$&$O_{8},$\\
 &$\vdots$&
\end{tabular}
\end{center}
Studying this multiplication table we can obtain some information about
the structure of $\mathrm{SLe(2)}$. For example,
there are some subalgebras important in the construction of
representations of $\mathrm{SLe(2)}$:
\begin{itemize}
 \item commutative negative grade subalgebra
    $A_{<0} = \left\langle O_1,E_2,E_3\right\rangle;$
 \item semisimple zero grade subalgebra
 $A_{0} = \left\langle E_4,E_5,E_6\right\rangle\simeq \mathrm{so(3)}
 \simeq\mathrm{sl(2)}\simeq\mathrm{sp(2)};$
 \item non-positive grade subalgebra
 $A_{\leq0} = A_{<0}+\hspace*{-10pt}\supset A_{0},$
 a semidirect sum of the semisimple algebra and the commutative ideal.
\end{itemize}
The results of computation of cohomology  $H^k_g(\mathrm{SLe(2)})$
are presented in Table \ref{tab:H1}.
The boxes of this table contain either three  numbers (with possible
indication of non-trivial cohomology class)
or right arrow.
The three numbers from top downwards are
\begin{itemize}
    \item $\dim C^k_g,$ the dimension of the whole space
    of $k$-cochains in grade $g$;
    \item the number of minimal subcomplexes
$C^{k-1}_{g,s}\stackrel{d^{k-1}_{g,s}} {\longrightarrow}
C^k_{g,s}\stackrel{d^k_{g,s}}{\longrightarrow}C^{k+1}_{g,s}$
constituting the whole subcomplex in accordance with formula (\ref{spart});
    \item $\max_{s\in S}\dim C^k_{g,s},$ maximum dimension of
    $(k,g)$-cochain subspaces
among all the minimal subcomplexes.
\end{itemize}
The right arrow $\rightarrow$ means that all subsequent boxes in
the row contain the same information. The reason for this is that all
the relations defining cocycles and coboundaries for $k$-cochains
coincide with those for $(k-1)$-cochains multiplied by the 1-cochain
$c(O_1)\equiv c(\theta\psi).$

In our computation we found four genuine cohomology classes,
i. e., generators of the cohomology ring, $\alpha,   \beta, \gamma, \delta.$
They are parenthesized in the table. Their multiplicative consequences are
underlined.
Notice that the cohomology ring contains nilpotents and zerodivisors:
there are arbitrary powers of the cocycle $\alpha = c(\theta\psi),$
but $\beta\alpha=0,$ $\gamma\alpha=0$ and $\delta\alpha^2=0.$

\noindent
\begin{table}[ht]
    \caption{$H^k_g(\mathrm{SLe(2)})$ for $(k,g) \in [1,\ldots,\infty)\times[-2k,\ldots,-2k+20]$}
    \label{tab:H1}
    \begin{center}
        \begin{tabular}{c|c|c|c|c|c|c|c|c|c||c}
\Tm{g+2k\backslash k}
    \Ti{1}\Ti{2}\Ti{3}\Ti{4}\Ti{5}\Ti{6}\Ti{7}\Ti{8}\Ti{9}\Tie{$k>9$}
\Ti{0} \Tc{1}{1}{1}{\left(\alpha\right)}
          \Tc{1}{1}{1}{\underline{\alpha^2}}
                \Tc{1}{1}{1}{\underline{\alpha^3}}
                      \Tc{1}{1}{1}{\underline{\alpha^4}}
                            \Tc{1}{1}{1}{\underline{\alpha^5}}
                                  \Tc{1}{1}{1}{\underline{\alpha^6}}
                                        \Tc{1}{1}{1}{\underline{\alpha^7}}
                                              \Tc{1}{1}{1}{\underline{\alpha^8}}
                                                    \Tc{1}{1}{1}{\underline{\alpha^9}}
                                                          \Tce{1}{1}{1}{\underline{\alpha^{k}}}\hline
\Ti{1} \Td{2}{2}{1}
          \Treprow
                \Treprow
                      \Treprow
                            \Treprow
                                  \Treprow
                                        \Treprow
                                              \Treprow
                                                    \Treprow
                                                          \Treprowe\hline
\Ti{2} \Td{3}{3}{1}
          \Td{4}{4}{1}
                \Treprow
                      \Treprow
                            \Treprow
                                  \Treprow
                                        \Treprow
                                              \Treprow
                                                    \Treprow
                                                          \Treprowe\hline
\Ti{3} \Td{6}{6}{1}
          \Td{12}{8}{3}
                \Treprow
                      \Treprow
                            \Treprow
                                  \Treprow
                                        \Treprow
                                              \Treprow
                                                    \Treprow
                                                          \Treprowe\hline
\Ti{4} \Td{8}{8}{1}
          \Tc{23}{13}{4}{\left(\beta\right)}
                \Td{26}{13}{4}
                      \Treprow
                            \Treprow
                                  \Treprow
                                        \Treprow
                                              \Treprow
                                                    \Treprow
                                                          \Treprowe\hline
\Ti{5} \Td{10}{10}{1}
          \Td{44}{16}{6}
                \Td{56}{18}{7}
                      \Treprow
                            \Treprow
                                  \Treprow
                                        \Treprow
                                              \Treprow
                                                    \Treprow
                                                          \Treprowe\hline
\Ti{6} \Td{12}{12}{1}
          \Td{73}{22}{8}
                \Td{118}{26}{15}
                      \Td{121}{26}{15}
                            \Treprow
                                  \Treprow
                                        \Treprow
                                              \Treprow
                                                    \Treprow
                                                          \Treprowe\hline
\Ti{7} \Td{14}{14}{1}
          \Td{116}{26}{10}
                \Td{226}{34}{23}
                      \Td{246}{34}{23}
                            \Treprow
                                  \Treprow
                                        \Treprow
                                              \Treprow
                                                    \Treprow
                                                          \Treprowe\hline
\Ti{8} \Td{16}{16}{1}
          \Td{171}{30}{14}
                \Td{414}{39}{36}
                      \Td{491}{42}{41}
                            \Td{492}{42}{41}
                                  \Treprow
                                        \Treprow
                                              \Treprow
                                                    \Treprow
                                                          \Treprowe\hline
\Ti{9} \Td{18}{18}{1}
          \Td{244}{34}{18}
                \Td{718}{48}{52}
                      \Td{952}{52}{71}
                            \Td{970}{52}{71}
                                  \Treprow
                                        \Treprow
                                              \Treprow
                                                    \Treprow
                                                          \Treprowe\hline
\Ti{10}\Td{20}{20}{1}
          \Td{333}{38}{23}
                \Td{1182}{54}{80}
                      \Td{1780}{65}{124}
                            \Tc{1867}{65}{124}{\left(\gamma\right)}
                                  \Td{1867}{65}{124}
                                        \Treprow
                                              \Treprow
                                                    \Treprow
                                                          \Treprowe\hline
\Ti{11}\Td{22}{22}{1}
          \Td{444}{42}{28}
                \Td{1870}{60}{119}
                      \Td{3204}{72}{197}
                            \Td{3528}{76}{197}
                                  \Td{3534}{76}{197}
                                        \Treprow
                                              \Treprow
                                                    \Treprow
                                                          \Treprowe\hline
\Ti{12}\Td{24}{24}{1}
          \Td{575}{46}{34}
                \Td{2858}{66}{176}
                      \Td{5584}{84}{311}
                            \Td{6546}{88}{358}
                                  \Td{6605}{88}{358}
                                        \Treprow
                                              \Treprow
                                                    \Treprow
                                                          \Treprowe\hline
\Ti{13}\Td{26}{26}{1}
          \Td{732}{50}{40}
                \Td{4224}{72}{241}
                      \Td{9398}{92}{489}
                            \Td{11878}{104}{606}
                                  \Td{12162}{104}{606}
                                        \Treprow
                                              \Treprow
                                                    \Treprow
                                                          \Treprowe\hline
\Ti{14}\Td{28}{28}{1}
          \Td{913}{54}{47}
                \Td{6082}{78}{330}
                      \Td{15343}{100}{787}
                            \Td{21073}{113}{1009}
                                  \Td{22102}{118}{1009}
                                        \Td{22119}{118}{1009}
                                              \Treprow
                                                    \Treprow
                                                          \Treprowe\hline
\Ti{15}\Td{30}{30}{1}
          \Td{1124}{58}{54}
                \Td{8552}{84}{434}
                      \Td{24348}{108}{1187}
                            \Td{36540}{130}{1578}
                                  \Td{39652}{134}{1598}
                                        \Td{39796}{134}{1598}
                                              \Treprow
                                                    \Treprow
                                                          \Treprowe\hline
\Ti{16}\Td{32}{32}{1}
          \Td{1363}{62}{62}
                \Td{11766}{90}{570}
                      \Td{37649}{116}{1776}
                            \Td{61884}{140}{2556}
                                  \Tc{70110}{153}{2802}{\left(\delta\right)}
                                        \Tc{70817}{153}{2802}{\underline{\delta\alpha}}
                                              \Td{70817}{153}{2802}
                                                    \Treprow
                                                          \Treprowe\hline
\Ti{17}\Td{34}{34}{1}
          \Td{1636}{66}{70}
                \Td{15892}{96}{721}
                      \Td{56848}{124}{2528}
                            \Td{102466}{150}{4127}
                                  \Td{122154}{164}{4610}
                                        \Td{124766}{170}{4610}
                                              \Td{124798}{170}{4610}
                                                    \Treprow
                                                          \Treprowe\hline
\Ti{18}\Td{36}{36}{1}
          \Td{1941}{70}{79}
                \Td{21114}{102}{912}
                      \Td{84034}{132}{3583}
                            \Td{165999}{160}{6594}
                                  \Td{209566}{184}{7510}
                                        \Td{217692}{188}{7510}
                                              \Td{217972}{188}{7510}
                                                    \Treprow
                                                          \Treprowe\hline
\Ti{19}\Td{38}{38}{1}
          \Td{2284}{74}{88}
                \Td{27622}{108}{1127}
                      \Td{121790}{140}{4900}
                            \Td{263372}{170}{10005}
                                  \Td{353762}{196}{11624}
                                        \Td{376010}{210}{11799}
                                              \Td{377422}{210}{11799}
                                                    \Treprow
                                                          \Treprowe\hline
\Ti{20}\Td{40}{40}{1}
          \Td{2663}{78}{98}
                \Td{35658}{114}{1392}
                      \Td{173394}{148}{6682}
                            \Td{409730}{180}{15095}
                                  \Td{587526}{208}{19154}
                                        \Td{642827}{234}{19996}
                                              \Td{648260}{246}{20438}
                                                    \Td{648308}{246}{20438}
                                                          \Treprowe\hline

        \end{tabular}
    \end{center}
\end{table}

\section{Concluding Remark}
When computing cohomology we start with the construction of the full set of $(k,g)$-monomials.
At the moment we do not see how to avoid this in deterministic algorithms.
To represent the set of monomials,
we need to allocate
$n=l\times\dim C^k_g$ elements of memory representing basis elements of algebra and module ($l=k$ for the trivial and $l=k+1$ for any non-trivial module).
In our implementation we represent basis elements by two-byte integers. For the last box in column 9 ($k = 9$) of Table \ref{tab:H1} we have $\dim C^k_g = 648308$ and the set of monomials occupies near 12MB. The dimensions grow very rapidly, so, in fact, we are working on the brink of abilities of 32bit architecture and, even theoretically, we can add only a few rows to Table \ref{tab:H1}. But, as to arithmetical difficulties,
we have some progress.

\section*{Acknowledgments}
I would like to thank H. Khudaverdian, D. Leites and I.
Shchepochkina for helpful comments on drafts of this text.

This work was partially supported by the
grant 01-01-00708 from the Russian Foundation for Basic Research and
grant 2339.2003.2 from the Russian Ministry of Industry, Science and
Technologies.

\end{document}